\documentclass[11pt]{article}
\usepackage{amsfonts}
\usepackage{amsfonts}
\usepackage{amsfonts}
\usepackage{mathrsfs}
\usepackage{bm}
\usepackage{amssymb,amsmath} % font used for R in Real numbers
 \allowdisplaybreaks

\catcode`!=11
\let\!int\int \def\int{\displaystyle\!int}
\let\!lim\lim \def\lim{\displaystyle\!lim}
\let\!sum\sum \def\sum{\displaystyle\!sum}
\let\!sup\sup \def\sup{\displaystyle\!sup}
\let\!inf\inf \def\inf{\displaystyle\!inf}
\let\!cap\cap \def\cap{\displaystyle\!cap}
\let\!max\max \def\max{\displaystyle\!max}
\let\!min\min \def\min{\displaystyle\!min}
\let\!frac\frac \def\frac{\displaystyle\!frac}
\catcode`!=12

\let\oldsection\section
\renewcommand\section{\setcounter{equation}{0}\oldsection}

\allowdisplaybreaks
\def\pf{\it{Proof.}\rm\quad}

\def\N{\mathbb{N}}\def\Z{\mathbb{Z}}

\newtheorem{thm}{Theorem}[section]
\newtheorem{lem}[thm]{Lemma}
\newtheorem{cor}[thm]{Corollary}
\newtheorem{exa}{Example}

\begin{document}
%%%%%%%%%%%%%%%%%%%% title %%%%%%%%%%%%%%%%%%%%%%%%%%%%%%%%%%%%%%%%%%%%%%%%
\title {\bf Explicit evaluation of quadratic Euler sums}
\author{
{Ce Xu\thanks{Corresponding author. Email: 15959259051@163.com},\quad Yingyue Yang,\quad Jianwen Zhang}\\[1mm]
\small School of Mathematical Sciences, Xiamen University\\
\small Xiamen
361005, P.R. China}

\date{}
\maketitle \noindent{\bf Abstract }  In this paper, we work out some explicit formulae for double nonlinear Euler sums involving harmonic numbers and alternating harmonic numbers. As applications of these formulae, we give new closed form
representations of several quadratic Euler sums through Riemann zeta function and linear sums. The given representations are new.
\\[2mm]
\\[2mm]

\section{Introduction}
In this paper, harmonic numbers, alternating harmonic numbers and their generalizations define as $$H_n=\sum\limits_{j=1}^n\frac {1}{j},\ \zeta_n(k)=\sum\limits_{j=1}^n\frac {1}{j^k} ,\ L_{n}(k)=\sum\limits_{j=1}^n\frac{(-1)^{j-1}}{j^k},\ 1\leq k \in \Z
.\eqno(1.1)$$
The classical linear Euler sum is defined by
$$S_{p,q}=\sum\limits_{n = 1}^\infty  {\frac{1}{{{n^q}}}} \sum\limits_{k = 1}^n {\frac{1}{{{k^p}}}}, \eqno(1.2)$$
where $p,q$ are positive integers, $q \geq 2$.
In 1742, Goldbach proposed to Euler the problem of expressing the $S_{p,q}$ in terms of values at positive integers of the Riemann zeta function $\zeta(s)$. The Riemann zeta function and  alternating Riemann zeta function are defined respectively by $$\zeta(s)=\sum\limits_{n = 1}^\infty {\frac {1}{n^{s}}},\Re(s)>1; \bar \zeta \left( s \right) = \sum\limits_{n = 1}^\infty  {\frac{{{{\left( { - 1} \right)}^{n - 1}}}}{{{n^s}}}} ,\;{\mathop{\Re}\nolimits} \left( s \right) \ge 1.$$
Euler showed this problem in the case $p = 1$ and gave a general formula for odd weight $p + q$  in 1775. He conjectured that the double linear
sums would be reducible to zeta values when p + q is odd, and even gave what he hoped to obtain the general
formula. In [3], D. Borwein, J.M. Borwein and R. Girgensohn proved conjecture and formula,
and in [2], D.H. Bailey, J.M. Borwein and R. Girgensohn conjectured that
the double linear sums when $p + q > 7,p + q$ is even, are not reducible. \\
Let $\pi  = \left( {{\pi _1}, \ldots ,{\pi _k}} \right)$ be a partition of integer $p$ and $p = {\pi _1} +  \cdots  + {\pi _k}$ with ${\pi _1} \le {\pi _2} \le  \cdots  \le {\pi _k}$. The classical double nonlinear Euler sum of index $\pi,q$ is defined as follows (see [12])
\[{S_{\pi ,q}} = \sum\limits_{n = 1}^\infty  {\frac{{{\zeta _n}\left( {{\pi _1}} \right){\zeta _n}\left( {{\pi _2}} \right) \cdots {\zeta _n}\left( {{\pi _k}} \right)}}{{{n^q}}}} ,\tag{1.3}\]
where the quantity ${\pi _1} +  \cdots  + {\pi _k} + q$ and $k$ called  the weight and the degree.\\
In general, we can also define the double Euler sums by the series
$$\sum\limits_{n = 1}^\infty  {\frac{{\prod\limits_{i = 1}^{{m_1}} {\zeta _n^{{q_{_i}}}\left( {{k_i}} \right)\prod\limits_{j = 1}^{{m_2}} {L_n^{{l_j}}\left( {{h_j}} \right)} } }}{{{n^p}}}} ,\ \sum\limits_{n = 1}^\infty  {\frac{{\prod\limits_{i = 1}^{{m_1}} {\zeta _n^{{q_{_i}}}\left( {{k_i}} \right)\prod\limits_{j = 1}^{{m_2}} {L_n^{{l_j}}\left( {{h_j}} \right)} {{\left( { - 1} \right)}^{n - 1}}} }}{{{n^p}}}}, \eqno(1.4)$$
where $p (p>1),m_{1},m_{2},q_{i},k_{i},h_{j},l_{j}$ are positive integers. If $\sum\limits_{i=1}^{m_{1}}(k_{i}q_{i})+\sum\limits_{j=1}^{m_{2}}(h_{j}l_{j})+p=C$($C$ is a positive integer), then we call the identity $C$th-order Euler sums.\\
In [12], Philippe Flajolet and Bruno Salvy gave explicit reductions to zeta values for all linear sums
\[\sum\limits_{n = 1}^\infty  {\frac{{{\zeta _n}\left( p \right)}}{{{n^q}}}} ,\;\sum\limits_{n = 1}^\infty  {\frac{{{L_n}\left( p \right)}}{{{n^q}}}} ,\;\sum\limits_{n = 1}^\infty  {\frac{{{\zeta _n}\left( p \right)}}{{{n^q}}}} {\left( { - 1} \right)^{n - 1}},\;\sum\limits_{n = 1}^\infty  {\frac{{{L_n}\left( p \right)}}{{{n^q}}}{{\left( { - 1} \right)}^{n - 1}}} \]
when $p+q$ is an odd weight. The evaluation of linear sums in
terms of values of Riemann zeta function and polylogarithm function at positive integers is known when $(p,q) = (1,3), (2,2)$, or $p + q$ is odd [6,12,18]. For instance, we have
$$\sum\limits_{n = 1}^\infty  {\frac{{{H_n}{}}}{{{n^3}}}}{\left( { - 1} \right)}^{n - 1}  =  - 2{\rm Li}{_4}\left( {\frac{1}{2}} \right) + \frac{{11{}}}{{4}} \zeta(4) + \frac{{{1}}}{{2}}\zeta(2){\ln ^2}2 - \frac{1}{{12}}{\ln ^4}2 - \frac{7}{4}\zeta \left( 3 \right)\ln 2,$$
$$\sum\limits_{n = 1}^\infty  {\frac{{{L_n}\left( 1 \right)}}{{{n^3}}}{{\left( { - 1} \right)}^{n - 1}} = \frac{3}{2}\zeta \left( 4 \right)}  + \frac{1}{2}\zeta \left( 2 \right){\ln ^2}2 - \frac{1}{{12}}{\ln ^4}2 - 2{\rm Li}{_4}\left( {\frac{1}{2}} \right).$$
The polylogarithm function defined for $\left| x \right| \le 1$ by
\[{\rm Li}{_p}\left( x \right) = \sum\limits_{n = 1}^\infty  {\frac{{{x^n}}}{{{n^p}}}}, \Re(p)>1 .\tag{1.5}\]
The relationship between the
values of the Riemann zeta function and Euler sums has been studied by many authors, for example see [2-4,6-17,19-20].\\
From [18], for a multi-index $S = \left( {{s_1},{s_2}, \cdots ,{s_k}} \right)\ \left( {{s_i} \in \N,\;{s_1} \ge 2} \right)$, the multiple zeta value (MZV for short) and the multiple zeta-star value (MZSV for short) are defined, respectively, by convergent series
\[\zeta \left( S \right) = \zeta \left( {{s_1},{s_2}, \cdots ,{s_k}} \right) = \sum\limits_{{n_1} > {n_2} >  \cdots  > {n_k} \ge 1} {\frac{1}{{n_1^{{s_1}}n_2^{{s_2}} \cdots n_k^{{s_k}}}}} ,\]
\[\zeta^\star\left( S \right) = \zeta \left( {{s_1},{s_2}, \cdots ,{s_k}} \right) = \sum\limits_{{n_1} \ge {n_2} \ge  \cdots  \ge {n_k} \ge 1} {\frac{1}{{n_1^{{s_1}}n_2^{{s_2}} \cdots n_k^{{s_k}}}}}.\]
We put a bar on top of $s_j\ (j=1,\cdots k)$ if there is a sign $(-1)^{n_j}$ appearing in the denominator on the right. For example,
\[{\zeta^\star}\left( {\bar 2,3} \right) = \sum\limits_{m \ge n \ge 1}^{} {\frac{{{{\left( { - 1} \right)}^{m }}}}{{{m^2}{n^3}}}}  = \sum\limits_{m = 1}^\infty  {\frac{{{{\left( { - 1} \right)}^{m }}}}{{{m^2}}}\sum\limits_{n = 1}^m {\frac{1}{{{n^3}}}} } .\]
By using the above notations, we have
\[{\zeta^\star}\left( {q,p} \right) = {S_{p,q}} = \sum\limits_{n = 1}^\infty  {\frac{{{\zeta _n}\left( p \right)}}{{{n^q}}}} ,\;{\zeta^\star}\left( {\bar q,p} \right) = -\sum\limits_{n = 1}^\infty  {\frac{{{\zeta _n}\left( p \right)}}{{{n^q}}}{{\left( { - 1} \right)}^{n - 1}}} ,\;\]
\[{\zeta^\star}\left( {q,\bar p} \right) = -\sum\limits_{n = 1}^\infty  {\frac{{{L_n}\left( p \right)}}{{{n^q}}}} ,\;{\zeta^\star}\left( {\bar q,\bar p} \right) = \sum\limits_{n = 1}^\infty  {\frac{{{L_n}\left( p \right)}}{{{n^q}}}{{\left( { - 1} \right)}^{n - 1}}} .\]
The main purpose of this paper is to evaluate some quadratic Euler sums which involving harmonic numbers and alternating harmonic numbers, either linearly or nonlinearly. In this paper, we will prove that all quadratic double sums
\[\sum\limits_{n = 1}^\infty  {\frac{{H_n^2}}{{{n^{2m}}}}} {\left( { - 1} \right)^{n - 1}},\sum\limits_{n = 1}^\infty  {\frac{{L_n^2\left( 1 \right)}}{{{n^{2m}}}}} {\left( { - 1} \right)^{n - 1}},\sum\limits_{n = 1}^\infty  {\frac{{{H_n}{L_n}\left( 1 \right)}}{{{n^{2m}}}}} {\left( { - 1} \right)^{n - 1}},\sum\limits_{n = 1}^\infty  {\frac{{{H_n}{L_n}\left( 1 \right)}}{{{n^{2m}}}}} ,\sum\limits_{n = 1}^\infty  {\frac{{L_n^2\left( 1 \right)}}{{{n^{2m}}}}} \]
are reducible to polynomials in zeta values and to linear sums, $m$ is a positive integer. Similarly, we have the following relations
\[\sum\limits_{n = 1}^\infty  {\frac{{H_n^2}}{{{n^{2m}}}}} {\left( { - 1} \right)^{n - 1}} = -2{\zeta^\star}\left( {\overline  {2m},1,1} \right) + {\zeta^\star}\left( {\overline  {2m},2} \right),\sum\limits_{n = 1}^\infty  {\frac{{L_n^2\left( 1 \right)}}{{{n^{2m}}}}} {\left( { - 1} \right)^{n - 1}} = -2{\zeta^\star}\left( {\overline { 2m},\overline 1,\overline 1} \right) + {\zeta^\star}\left( {\overline  {2m},2} \right),\]
\[\sum\limits_{n = 1}^\infty  {\frac{{L_n^2\left( 1 \right)}}{{{n^{2m}}}}}  = 2{\zeta^\star}\left( {2m,\overline1,\overline1} \right) - {\zeta^\star}\left( {2m,2} \right),\;\sum\limits_{n = 1}^\infty  {\frac{{\zeta _n^2\left( p \right)}}{{{n^{2m}}}}}  = 2{\zeta^\star}\left( {2m,p,p} \right) - {\zeta^\star}\left( {2m,2p} \right),\]
\[\sum\limits_{n = 1}^\infty  {\frac{{{H_n}{L_n}\left( 1 \right)}}{{{n^{2m}}}}}  = -{\zeta^\star}\left( {2m,\overline1,1} \right) - {\zeta^\star}\left( {2m,1,\overline1} \right) + {\zeta^\star}\left( {2m,\overline2} \right),\]
\[\sum\limits_{n = 1}^\infty  {\frac{{{H_n}{L_n}\left( 1 \right)}}{{{n^{2m}}}}} {\left( { - 1} \right)^{n - 1}} = {\zeta^\star}\left( {\overline{2m},\overline1,1} \right) + {\zeta^\star}\left( {\overline{2m},1,\overline1} \right) - {\zeta^\star}\left( {\overline{2m},\overline2} \right).\]
The following lemma will be useful in the development of the main theorems.
\begin{lem}\label{lem 2.1}
Let $n$ be a positive integers. We have the following relations
\[\sum\limits_{k = 1}^{n - 1} {\frac{{{L_k}\left( 1 \right)}}{{n - k}}{{\left( { - 1} \right)}^{k - 1}}}  = 2{\left( { - 1} \right)^{n - 1}}\sum\limits_{k = 1}^n {\frac{{{H_k}}}{k}{{\left( { - 1} \right)}^{k - 1}}}  - 2{\left( { - 1} \right)^{n - 1}}{L_n}\left( 2 \right),\tag{1.6}\]
\[\sum\limits_{m = n + 1}^\infty  {\frac{{{L_m}\left( 1 \right)}}{{m - n}}{{\left( { - 1} \right)}^{m - n}}}  = \sum\limits_{k = 1}^n {\frac{{{H_k}}}{k}{{\left( { - 1} \right)}^{k - 1}}}  - \sum\limits_{j = 1}^\infty  {\frac{{{L_j}\left( 1 \right)}}{j}} {\left( { - 1} \right)^{j - 1}},\tag{1.7}\]
\[{\left( { - 1} \right)^{n - 1}}\sum\limits_{k = 1}^{n - 1} {\frac{{{L_k}\left( 1 \right)}}{{n - k}}}  = \sum\limits_{k = 1}^{n - 1} {\frac{{{H_k}}}{{n - k}}{{\left( { - 1} \right)}^k}} ,\tag{1.8}\]
\[\sum\limits_{m = n + 1}^\infty  {\frac{{{H_m}}}{{m - n}}{{\left( { - 1} \right)}^{m - n}}}  = \frac{{L_n^2\left( 1 \right) + {\zeta _n}\left( 2 \right)}}{2} - \ln 2\left( {{H_n} + {L_n}\left( 1 \right)} \right) - \sum\limits_{n = 1}^\infty  {\frac{{{H_j}}}{j}} {\left( { - 1} \right)^{j - 1}}.\tag{1.9}\]
\end{lem}
where $\sum\limits_{n{\rm{ = }}1}^\infty  {\frac{{{H_n}{{\left( { - 1} \right)}^{n - 1}}}}{n}}  = \frac{{\zeta \left( 2 \right) - {{\ln }^2}2}}{2},\sum\limits_{n{\rm{ = }}1}^\infty  {\frac{{{L_n}(1){{\left( { - 1} \right)}^{n - 1}}}}{n}}  = \frac{{\zeta \left( 2 \right) + {{\ln }^2}2}}{2}.$\\
\pf To prove (1.6), using Cauchy product formula, we have
\[\frac{{\ln \left( {1 + x} \right)}}{{1 - x}} = \left( {\sum\limits_{n = 1}^\infty  {\frac{{{{\left( { - 1} \right)}^{n - 1}}}}{n}{x^n}} } \right)\left( {\sum\limits_{n = 1}^\infty  {{x^{n - 1}}} } \right) = \sum\limits_{n = 1}^\infty  {{L_n}\left( 1 \right){x^n}}.\tag{1.10}\]
Furthermore, by (1.10) and Cauchy product formula, we can obtain
\[\frac{{{{\ln }^2}\left( {1 + x} \right)}}{{1 - x}} = \sum\limits_{n = 1}^\infty  {{{\left( { - 1} \right)}^{n + 1}}\left( {\sum\limits_{k = 1}^n {\frac{{{L_k}\left( 1 \right)}}{{n - k + 1}}{{\left( { - 1} \right)}^{k - 1}}} } \right){x^{n + 1}}}.\tag{1.11}\]
On the other hand, from [14], we get
\[{\ln ^2}\left( {1 + x} \right) = 2\left\{ {\sum\limits_{n = 1}^\infty  {\frac{{{{\left( { - 1} \right)}^{n - 1}}}}{{{n^2}}}{x^n}}  - \sum\limits_{n = 1}^\infty  {\frac{{{H_n}}}{n}{{\left( { - 1} \right)}^{n - 1}}{x^n}} } \right\}.\tag{1.12}\]
Similarly, applying (1.12), then have
\[\frac{{{{\ln }^2}\left( {1 + x} \right)}}{{1 - x}} = 2\sum\limits_{n = 1}^\infty  {\sum\limits_{k = 1}^n {\left\{ {{L_n}\left( 2 \right) - \frac{{{H_k}}}{k}{{\left( { - 1} \right)}^{k - 1}}} \right\}} } {x^n}.\tag{1.13}\]
Comparing (1.11) and (1.13), we can obtain (1.6). Now, we prove (1.7). We see that we may rewrite the left side of (1.7) as
\begin{align*}
\sum\limits_{m = n + 1}^\infty  {\frac{{{L_m}\left( 1 \right)}}{{m - n}}{{\left( { - 1} \right)}^{m - n}}}
& =\sum\limits_{j = 1}^\infty  {\frac{{{L_{n + j}}\left( 1 \right)}}{j}{{\left( { - 1} \right)}^j}} \\
& =\sum\limits_{j = 1}^\infty  {\frac{{{{\left( { - 1} \right)}^j}}}{j}\left\{ {\sum\limits_{k = 1}^j {\frac{{{{\left( { - 1} \right)}^{k - 1}}}}{k}}  + \sum\limits_{k = j + 1}^{n + j} {\frac{{{{\left( { - 1} \right)}^{k - 1}}}}{k}} } \right\}}
\\
& =\sum\limits_{j = 1}^\infty  {\frac{{{L_j}\left( 1 \right)}}{j}{{\left( { - 1} \right)}^j}}  + \sum\limits_{k = 1}^n {\frac{{{{\left( { - 1} \right)}^{k - 1}}}}{k}\sum\limits_{j = 1}^\infty  {\left( {\frac{1}{j} - \frac{1}{{j + k}}} \right)} }
\\
& =\sum\limits_{k = 1}^n {\frac{{{H_k}}}{k}{{\left( { - 1} \right)}^{k - 1}}}  - \sum\limits_{j = 1}^\infty  {\frac{{{L_j}\left( 1 \right)}}{j}} {\left( { - 1} \right)^{j - 1}} .  \tag{1.14}
\end{align*}
Similarly to the proof of (1.6) and (1.7), we can prove (1.8) and (1.9).
\section{Main Theorems and Proof}
In this section, we will establish some explicit relationships which
involve alternating quadratic Euler sums and linear sums.
\begin{thm}\label{lem 2.1}
For integer $m\geq 1$, we have
\begin{align*}
&\sum\limits_{j = 1}^{m - 1} {{{\left( { - 1} \right)}^{j - 1}}\bar \zeta \left( {m - j + 1} \right){\zeta^\star}\left( {j+1,1} \right)}  + {\left( { - 1} \right)^{m - 1}}\ln 2{\zeta^\star}\left( {m+1,1} \right)  - {\zeta^\star}\left( {\overline {m+2},1} \right) \\
&\quad - {\left( { - 1} \right)^{m - 1}}\ln 2{\zeta^\star}\left( {\overline {m+1},1} \right)  - {\left( { - 1} \right)^{m - 1}}\sum\limits_{n = 1}^\infty  {\frac{{{{\left( { - 1} \right)}^{n - 1}}{H_n}{L_n}\left( 1 \right)}}{{{n^{m + 1}}}}} \\
&\quad = \frac{1}{2}\sum\limits_{n = 1}^\infty  {\frac{{{{\left( { - 1} \right)}^{n - 1}}H_n^2}}{{{n^{m + 1}}}}}  -\frac{1}{2}{\zeta^\star}\left( {\overline {m+1},2} \right)  + \zeta \left( 2 \right)\bar \zeta \left( {m + 1} \right).\tag{2.1}
\end{align*}
\end{thm}
\pf
Let$${I_m} = \sum\limits_{n = 1}^\infty  {\frac{{{{\left( { - 1} \right)}^{n - 1}}}}{{{n^m}\left( {n + k} \right)}}}.$$
Using the definition of $I_m$, we can easily verify that
\[I_m = \sum\limits_{j = 1}^{m - 1} {\frac{{{{\left( { - 1} \right)}^{j - 1}}}}{{{k^j}}}\bar \zeta \left( {m - j + 1} \right)} \; + \frac{{{{\left( { - 1} \right)}^{m - 1}}}}{{{k^m}}}\ln 2 + \frac{{{{\left( { - 1} \right)}^{m + k}}}}{{{k^m}}}\ln 2 - \frac{{{{\left( { - 1} \right)}^{m + k}}}}{{{k^m}}}{L_k}\left( 1 \right).\tag{2.2}\]
On the other hand, we have
$$\sum\limits_{k = 1}^\infty  {\frac{{{H_k}}}{{k\left( {n + k} \right)}}}  = \frac{1}{n}\left( {\frac{1}{2}H_n^2 + \frac{1}{2}{\zeta _n}\left( 2 \right) + \zeta \left( 2 \right) - \frac{{{H_n}}}{n}} \right).\eqno(2.3)$$
Multiplying (2.3) by $\frac{{{{\left( { - 1} \right)}^{n - 1}}}}{{{n^m}}}$ and summing with respect to $n$. The result is
\[\sum\limits_{n = 1}^\infty  {\sum\limits_{k = 1}^\infty  {\frac{{{{\left( { - 1} \right)}^{n - 1}}{H_k}}}{{k{n^m}\left( {n + k} \right)}}} }  = \frac{1}{2}\sum\limits_{n = 1}^\infty  {\frac{{{{\left( { - 1} \right)}^{n - 1}}H_n^2}}{{{n^{m + 1}}}}}  - \frac{1}{2}{\zeta^\star}\left( {\overline {m+1},2} \right) \; + \zeta \left( 2 \right)\bar \zeta \left( {m + 1} \right) + {\zeta^\star}\left( {\overline {m+2},1} \right) .\tag{2.4}\]
Substituting (2.2) into (2.4) yields desired result.
\begin{thm}\label{lem 2.2}\ \ Let $s \ge 1$ be integer. Then the following identity holds
\begin{align*}
\frac{1}{2}\sum\limits_{n = 1}^\infty  {\frac{{H_n^2}}{{{n^s}}}{{\left( { - 1} \right)}^{n - 1}}}  - \sum\limits_{n = 1}^\infty  {\frac{{{H_n}{L_n}\left( 1 \right)}}{{{n^s}}}{{\left( { - 1} \right)}^{n - 1}}}
 &=-\frac{3}{2}{\zeta^\star}\left( {\overline {s},2} \right)  - s{\zeta^\star}\left( {\overline {s+1},1} \right)  - {\zeta^\star}\left( { {s+1},1} \right)
\nonumber \\
& \quad- \sum\limits_{j = 2}^s { { {\zeta^\star}\left( { {j},1} \right)} \bar \zeta \left( {s + 1 - j} \right)} + \ln 2{\zeta^\star}\left( {\overline {s},1} \right)
,\tag{2.5}
\end{align*}
\[\frac{3}{2}\sum\limits_{n = 1}^\infty  {\frac{{H_n^2 - {\zeta _n}\left( 2 \right)}}{{{n^s}}}}  = \left( {s + 1} \right){\zeta^\star}\left( { {s+1},1} \right)  - \sum\limits_{j = 2}^{s - 1} {{\zeta^\star}\left( { {j},1} \right)\zeta \left( {s + 1 - j} \right)},\;2 \le s \in Z.
\tag{2.6}\]
\end{thm}
\pf For real $-1\leq x<1$ and integers $s\ (s\geq 2)$ and $t\geq 0$, consider function
\begin{align*}
{T_{s,t}}\left( {x,y} \right)
 &=\sum\limits_{\scriptstyle n,m = 1 \hfill \atop
  \scriptstyle n \ne m \hfill}^\infty  {\frac{{{H_n}{x^n}{y^m}}}{{{n^s}{m^t}\left( {m - n} \right)}}}  = \sum\limits_{\scriptstyle n,m = 1 \hfill \atop
  \scriptstyle n \ne m \hfill}^\infty  {\frac{{{H_n}{x^n}{y^m}\left( {m - n + n} \right)}}{{{n^s}{m^{t + 1}}\left( {m - n} \right)}}}
\nonumber \\
& = \sum\limits_{\scriptstyle n,m = 1 \hfill \atop
  \scriptstyle n \ne m \hfill}^\infty  {\frac{{{H_n}{x^n}{y^m}}}{{{n^s}{m^{t + 1}}}}}  + \sum\limits_{\scriptstyle n,m = 1 \hfill \atop
  \scriptstyle n \ne m \hfill}^\infty  {\frac{{{H_n}{x^n}{y^m}}}{{{n^{s - 1}}{m^{t + 1}}\left( {m - n} \right)}}}
\nonumber \\
& =  \sum\limits_{n = 1}^\infty  {\frac{{{H_n}{x^n}}}{{{n^s}}}} \left( {\sum\limits_{m = 1}^\infty  {\frac{{{y^m}}}{{{m^{t + 1}}}}}  - \frac{{{y^n}}}{{{n^{t + 1}}}}} \right) + {T_{s - 1,t + 1}}\left( {x,y} \right)
\nonumber \\
& = \left( {\sum\limits_{n = 1}^\infty  {\frac{{{H_n}{x^n}}}{{{n^s}}}} } \right)L{i_{t + 1}}\left( y \right) - \left( {\sum\limits_{n = 1}^\infty  {\frac{{{H_n}{x^n}{y^n}}}{{{n^{s + t + 1}}}}} } \right) + {T_{s - 1,t + 1}}\left( {x,y} \right)
.\tag{2.7}
\end{align*}
Define
\[H_s^{\left( k \right)}\left( x \right) = \;\;\sum\limits_{n = 1}^\infty  {\frac{{H_n^k{x^n}}}{{{n^s}}}} ,\; - 1 \le x < 1,\;0 \le k \in \Z, \Re(s)\geq1.\]
Obviously, $H_s^{\left( 0 \right)}\left( x \right) = L{i_s}\left( x \right)$. Hence, we can rewrite (2.7) as
\[{T_{s,t}}\left( {x,y} \right) = H_s^{\left( 1 \right)}\left( x \right)H_{t + 1}^{\left( 0 \right)}\left( y \right) - H_{s + t + 1}^{\left( 1 \right)}\left( {xy} \right) + {T_{s - 1,t + 1}}\left( {x,y} \right).\tag{2.8}\]
Telescoping this gives
\[{T_{s,t}}\left( {x,y} \right) = \sum\limits_{j = 1}^s {H_j^{\left( 1 \right)}\left( x \right)H_{s + t + 1 - j}^{\left( 0 \right)}\left( y \right)}  - sH_{s + t + 1}^{\left( 1 \right)}\left( {xy} \right) + {T_{0,s + t}}\left( {x,y} \right).\tag{2.9}\]
Taking $t=0$ in (2.9) yield
\[{T_{s,0}}\left( {x,y} \right) = \sum\limits_{j = 1}^s {H_j^{\left( 1 \right)}\left( x \right)H_{s + 1 - j}^{\left( 0 \right)}\left( y \right)}  - sH_{s + 1}^{\left( 1 \right)}\left( {xy} \right) + {T_{0,s}}\left( {x,y} \right).\tag{2.10}\]
From the definition of ${T_{s,t}}\left( {x,y} \right)$, we conclude that
\begin{align*}
{T_{s,0}}\left( {x,y} \right)
 & = \sum\limits_{\scriptstyle n,m = 1 \hfill \atop
  \scriptstyle n \ne m \hfill}^\infty  {\frac{{{H_n}{x^n}{y^m}}}{{{n^s}\left( {m - n} \right)}}}  = \sum\limits_{n = 1}^\infty  {\frac{{{H_n}{x^n}}}{{{n^s}}}} \sum\limits_{\scriptstyle m = 1 \hfill \atop
  \scriptstyle m \ne n \hfill}^\infty  {\frac{{{y^m}}}{{m - n}}}
\nonumber \\
&  = \sum\limits_{n = 1}^\infty  {\frac{{{H_n}{x^n}{y^n}}}{{{n^s}}}} \sum\limits_{m = n + 1}^\infty  {\frac{{{y^{m - n}}}}{{m - n}}}  - \sum\limits_{n = 1}^\infty  {\frac{{{H_n}{x^n}}}{{{n^s}}}\sum\limits_{m = 1}^{n - 1} {\frac{{{y^m}}}{{n - m}}} }
\nonumber \\
&  = H_s^{\left( 1 \right)}\left( {xy} \right)H_1^{\left( 0 \right)}\left( y \right) - \sum\limits_{n = 1}^\infty  {\frac{{{H_n}{x^n}}}{{{n^s}}}\left( {\sum\limits_{j = 1}^{n - 1} {\frac{{{y^{n - j}}}}{j}} } \right)}
.\tag{2.11}
\end{align*}
\begin{align*}
{T_{0,s}}\left( {x,y} \right)
 &  = \sum\limits_{\scriptstyle n,m = 1 \hfill \atop
  \scriptstyle n \ne m \hfill}^\infty  {\frac{{{H_n}{x^n}{y^m}}}{{{m^s}\left( {m - n} \right)}}}  = \sum\limits_{m = 1}^\infty  {\frac{{{y^m}}}{{{m^s}}}\sum\limits_{\scriptstyle n = 1 \hfill \atop
  \scriptstyle n \ne m \hfill}^\infty  {\frac{{{H_n}{x^n}}}{{m - n}}} }
\nonumber \\
&   =  - \sum\limits_{m = 1}^\infty  {\frac{{{x^m}{y^m}}}{{{m^s}}}\sum\limits_{n = m + 1}^\infty  {\frac{{{H_n}{x^{n - m}}}}{{n - m}}} }  +  \sum\limits_{m = 1}^\infty  {\frac{{{y^m}}}{{{m^s}}}\sum\limits_{n = 1}^{m - 1} {\frac{{{H_n}{x^n}}}{{m - n}}} }
.\tag{2.12}
\end{align*}
Noting that
\begin{align*}
\sum\limits_{n = m + 1}^\infty  {\frac{{{H_n}{x^{n - m}}}}{{n - m}}}
 & = \sum\limits_{j = 1}^\infty  {\frac{{{H_{n + j}}}}{j}{x^j}}  = \sum\limits_{j = 1}^\infty  {\frac{{{x^j}}}{j}} \left( {\sum\limits_{k = 1}^j {\frac{1}{k}}  + \sum\limits_{k = j + 1}^{n + j} {\frac{1}{k}} } \right)
\nonumber \\
&  = \sum\limits_{j = 1}^\infty  {\frac{{{H_j}}}{j}{x^j}}  + \sum\limits_{j = 1}^\infty  {\frac{{{x^j}}}{j}} \sum\limits_{k = 1}^n {\frac{1}{{k + j}}}
\nonumber \\
&  = H_1^{\left( 1 \right)}\left( x \right)+ \sum\limits_{k = 1}^n {\frac{1}{k}} \sum\limits_{j = 1}^\infty  {{x^j}\left( {\frac{1}{j} - \frac{1}{{k + j}}} \right)}
,\tag{2.13}
\end{align*}
\[\sum\limits_{j = 1}^{n - 1} {\frac{{{H_j}}}{{n - j}}}  = H_n^2 - {\zeta _n}\left( 2 \right),\;\sum\limits_{k = 1}^n {\frac{{{H_k}}}{k}}  = \frac{{H_n^2 + {\zeta _n}\left( 2 \right)}}{2}.\tag{2.14}\]
\[\mathop {\lim }\limits_{x \to 1,y \to 1} \left\{ {H_s^{\left( 1 \right)}\left( {xy} \right)H_1^{\left( 0 \right)}\left( y \right) - H_s^{\left( 1 \right)}\left( x \right)H_1^{\left( 0 \right)}\left( y \right)} \right\} = 0,\tag{2.15}\]
\[\mathop {\lim }\limits_{x \to 1,y \to 1} \left\{ {H_s^{\left( 0 \right)}\left( {xy} \right)H_1^{\left( 1 \right)}\left( x \right) - H_s^{\left( 0 \right)}\left( y \right)H_1^{\left( 1 \right)}\left( x \right)} \right\} = 0,\tag{2.16}\]
\[\mathop {\lim }\limits_{x \to 1,y \to  - 1} \left\{ {H_s^{\left( 0 \right)}\left( {xy} \right)H_1^{\left( 1 \right)}\left( x \right) - H_s^{\left( 0 \right)}\left( y \right)H_1^{\left( 1 \right)}\left( x \right)} \right\} = 0.\tag{2.17}\]
Combining (2.10)-(2.17) and taking $(x,y)\rightarrow (1,1), (1,-1)$, we can obtain (2.5) and (2.6). If taking $(x,y)\rightarrow (1,-1)$ in (2.10)-(2.12) and combining (1.9), we can give the following Corollary.
\begin{cor}
Let $s\geq 2$ be integer, we have
\begin{align*}
\sum\limits_{n = 1}^\infty  {\frac{{\sum\limits_{k = 1}^{n - 1} {\frac{{{H_k}}}{{n - k}}{{\left( { - 1} \right)}^k}} }}{{{n^s}}}}
 &=\sum\limits_{n = 1}^\infty  {\frac{{H_n^2}}{{{n^s}}}{{\left( { - 1} \right)}^{n - 1}}}  + \frac{1}{2}\left( {\zeta \left( 2 \right) - {{\ln }^2}2} \right)\bar \zeta \left( s \right) + \ln 2\left({\zeta^\star}\left( {\overline s,\overline1} \right)-{\zeta^\star}\left( {\overline s,1} \right) \right)
  \nonumber \\
           &\quad - \sum\limits_{j = 1}^{s - 1} {{\zeta^\star}\left( {\overline j,1} \right)\zeta \left( {s + 1 - j} \right)} + \left( {s + 1} \right){\zeta^\star}\left( {\overline {s+1},1} \right)
  \nonumber \\
           & \quad  - \frac{1}{2}\sum\limits_{n = 1}^\infty  {\frac{{L_n^2\left( 1 \right) }}{{{n^s}}}{{\left( { - 1} \right)}^{n - 1}}} +\frac{1}{2}{\zeta^\star}\left( {\overline {s},2} \right) .\tag{2.18}
\end{align*}
\end{cor}
In the same way, we can obtain the following Theorem.
\begin{thm} Let $s\geq 1$ be integer, then have
\begin{align*}
{\zeta^\star}\left( {s,\bar 1,1} \right)&=-\sum\limits_{n = 1}^\infty  {\frac{{\sum\limits_{k = 1}^n {\frac{{{H_k}}}{k}{{\left( { - 1} \right)}^{k - 1}}} }}{{{n^s}}}}\\
 &=s{\zeta^\star}\left( {{s+1},\overline 1} \right) - \sum\limits_{n = 1}^\infty  {\frac{{L_n^2\left( 1 \right)}}{{{n^s}}}} + 2{\zeta^\star}\left( {{s},\overline 2} \right) - \ln 2{\zeta^\star}\left( {{s},\overline 1} \right)
  \nonumber \\
           &\quad + {\zeta^\star}\left( {\overline{s+1},\overline 1} \right) + \zeta \left( s \right){\zeta^\star}\left( {\overline{1},\overline 1} \right)+ \sum\limits_{j = 1}^s {{\zeta^\star}\left( {\overline{j},\overline 1} \right)\bar \zeta \left( {s + 1 - j} \right)},\ (s\geq 2)  .\tag{2.19}
\end{align*}
\begin{align*}
{\zeta^\star}\left( {\overline{s},\overline 1,1} \right)&=\sum\limits_{n = 1}^\infty  {\frac{{\sum\limits_{k = 1}^n {\frac{{{H_k}}}{k}{{\left( { - 1} \right)}^{k - 1}}} }}{{{n^s}}}{{\left( { - 1} \right)}^{n - 1}}}\\
 & =\left( {s + 1} \right){\zeta^\star}\left( {\overline{s+1},\overline 1} \right)  + 2{\zeta^\star}\left( {\overline{s},\overline 2} \right)- \bar \zeta \left( s \right){\zeta^\star}\left( {\overline{1},\overline 1} \right)  - \sum\limits_{n = 1}^\infty  {\frac{{{H_n}{L_n}\left( 1 \right)}}{{{n^s}}}{{\left( { - 1} \right)}^{n - 1}}}
\nonumber \\
&\quad  - \sum\limits_{j = 1}^{s - 1} {{\zeta^\star}\left( {\overline{j},\overline 1} \right) \zeta \left( {s + 1 - j} \right)} .\tag{2.20}
\end{align*}
\begin{align*}
\sum\limits_{n = 1}^\infty  {\frac{{{{\left( { - 1} \right)}^{n - 1}}\sum\limits_{k = 1}^{n - 1} {\frac{{{L_k}\left( 1 \right)}}{{n - k}}} }}{{{n^s}}}}
 & =s{\zeta^\star}\left( {\overline{s+1},\overline 1} \right)   + \sum\limits_{n = 1}^\infty  {\frac{{L_n^2\left( 1 \right)}}{{{n^s}}}{{\left( { - 1} \right)}^{n - 1}}} + {\zeta^\star}\left( {{s+1},\overline 1} \right)
\nonumber \\
&\quad + \sum\limits_{j = 2}^s {{\zeta^\star}\left( {{j},\overline 1} \right) \bar \zeta \left( {s + 1 - j} \right)}  - \ln 2{\zeta^\star}\left( {\overline{s},\overline 1} \right)
\nonumber \\
&\quad - \ln 2\left({\zeta^\star}\left( {\overline{s},\overline 1} \right)-{\zeta^\star}\left( {\overline{s}, 1} \right)  \right)  + {\zeta^\star}\left( {\overline{s+1},1,\overline 1} \right) .\tag{2.21}
\end{align*}
\end{thm}
\pf Similarly to the proof of Theorem 2.2, considering function
\[{T_{s,t}}\left( {x,y} \right) = \sum\limits_{\scriptstyle n,m = 1 \hfill \atop
  \scriptstyle n \ne m \hfill}^\infty  {\frac{{{L_n}\left( 1 \right){x^n}{y^m}}}{{{n^s}{m^t}\left( {m - n} \right)}}} .\tag{2.22}\]
We have the following identity
\begin{align*}
&\sum\limits_{n = 1}^\infty  {\frac{{{L_n}\left( 1 \right)}}{{{n^s}}}{x^n}{y^n}} \sum\limits_{n = 1}^\infty  {\frac{{{y^n}}}{n}}  - \sum\limits_{n = 1}^\infty  {\frac{{{L_n}\left( 1 \right)}}{{{n^s}}}{x^n}} \sum\limits_{j = 1}^{n - 1} {\frac{{{y^{n - j}}}}{j}}\\
&\quad + \sum\limits_{n = 1}^\infty  {\frac{{{x^n}{y^n}}}{{{n^s}}}} \sum\limits_{m = n + 1}^\infty  {\frac{{{L_m}\left( 1 \right){x^{m - n}}}}{{m - n}}}  - \sum\limits_{n = 1}^\infty  {\frac{{{y^n}}}{{{n^s}}}} \sum\limits_{m = 1}^{n - 1} {\frac{{{L_m}\left( 1 \right){x^m}}}{{n - m}}}  \\
&\quad = \sum\limits_{j = 1}^s {\left( {\sum\limits_{n = 1}^\infty  {\frac{{{L_n}\left( 1 \right)}}{{{n^j}}}{x^n}} } \right)\left( {\sum\limits_{n = 1}^\infty  {\frac{{{y^n}}}{{{n^{s + 1 - j}}}}} } \right)}  - s\sum\limits_{n = 1}^\infty  {\frac{{{L_n}\left( 1 \right)}}{{{n^{s + 1}}}}{x^n}{y^n}} .\tag{2.23}
\end{align*}
Noting that
\begin{align*}
\sum\limits_{m = n + 1}^\infty  {\frac{{{L_m}\left( 1 \right)}}{{m - n}}} {x^{m - n}}
 &= \sum\limits_{j = 1}^\infty  {\frac{{{L_{n + j}}\left( 1 \right)}}{j}{x^j}}  = \sum\limits_{j = 1}^\infty  {\frac{{{x^j}}}{j}} \left\{ {\sum\limits_{k = 1}^j {\frac{{{{\left( { - 1} \right)}^{k - 1}}}}{k}}  + \sum\limits_{k = j + 1}^{n + j} {\frac{{{{\left( { - 1} \right)}^{k - 1}}}}{k}} } \right\}
\nonumber \\
&  = \sum\limits_{j = 1}^\infty  {\frac{{{L_j}\left( 1 \right)}}{j}{x^j}}  + \sum\limits_{k = 1}^n {{{\left( { - 1} \right)}^{k - 1}}\sum\limits_{j = 1}^\infty  {\frac{{{{\left( { - x} \right)}^j}}}{{j\left( {k + j} \right)}}} }
\nonumber \\
&  =\sum\limits_{j = 1}^\infty  {\frac{{{L_j}\left( 1 \right)}}{j}{x^j}}  + \sum\limits_{k = 1}^n {\frac{{{{\left( { - 1} \right)}^{k - 1}}}}{k}\sum\limits_{j = 1}^\infty  {{{\left( { - x} \right)}^j}\left( {\frac{1}{j} - \frac{1}{{k + j}}} \right)} }
\tag{2.24}
\end{align*}
and
\[\mathop {\lim }\limits_{x \to 1} \sum\limits_{k = 1}^n {\frac{{{{\left( { - 1} \right)}^{k - 1}}}}{k}\sum\limits_{j = 1}^\infty  {{{\left( { - x} \right)}^j}\left( {\frac{1}{j} - \frac{1}{{k + j}}} \right)} }  = \sum\limits_{k = 1}^n {\frac{{{L_k}\left( 1 \right)}}{k}}  - \ln 2\left( {{H_n} + {L_n}\left( 1 \right)} \right).\tag{2.25}\]
Taking $(x,y)\rightarrow(-1,-1),(-1,1),(1,-1)$ and combining (1.6)(1.7)(2.23)(2.24), we can obtain (2.19)-(2.21). Next, we give two explicit relationships which involve alternating double Euler sums and triple linear sums.
\begin{thm}
For $p\geq 1$ and $p\in \Z$, we have
\begin{align*}
{\zeta^\star}\left( {\overline{p+1},1,\overline 1} \right) &=\sum\limits_{n = 1}^\infty  {\frac{{{{\left( { - 1} \right)}^{n - 1}}}}{{{n^{p + 1}}}}\left( {\sum\limits_{k = 1}^n {\frac{{{L_k}\left( 1 \right)}}{k}} } \right)}\\
& =\sum\limits_{i = 1}^{p - 1} {{{\left( { - 1} \right)}^{i - 1}}\zeta \left( {p + 1 - i} \right)\bar \zeta \left( {i + 2} \right)}  - {\left( { - 1} \right)^{p - 1}}{\zeta^\star}\left( {\overline{p+2},1} \right)  \\
&\quad  + \ln 2\left({\zeta^\star}\left( {\overline{p+1},\overline 1} \right)-{\zeta^\star}\left( {\overline{p+1},1} \right)\right)- \frac{1}{2}{\ln ^2}2\left\{ {\zeta \left( {p + 1} \right) + \bar \zeta \left( {p + 1} \right)} \right\}   \\
&\quad + \sum\limits_{i = 1}^{p - 1} {{{\left( { - 1} \right)}^{i - 1}}\zeta \left( {p + 1 - i} \right){\zeta^\star}\left( {\overline{i+1},1} \right)}
  \\
&\quad  - {\left( { - 1} \right)^{p - 1}}\sum\limits_{n = 1}^\infty  {\frac{{H_n^2}}{{{n^{p + 1}}}}{{\left( { - 1} \right)}^{n - 1}}} .\tag{2.26}
\end{align*}
\begin{align*}
-{\zeta^\star}\left( {{p+1},1,\overline 1} \right)&=\sum\limits_{n = 1}^\infty  {\frac{1}{{{n^{p + 1}}}}\left( {\sum\limits_{k = 1}^n {\frac{{{L_k}\left( 1 \right)}}{k}} } \right)}\\
& =\ln 2\left({\zeta^\star}\left( {{p+1},1} \right)-{\zeta^\star}\left( {{p+1},\overline 1} \right)\right) + \sum\limits_{i = 1}^{p - 1} {{{\left( { - 1} \right)}^{i - 1}}\bar \zeta \left( {p + 1 - i} \right){\zeta^\star}\left( {\overline{i+1}, 1} \right) } \\
&\quad + \ln 2{\left( { - 1} \right)^{p - 1}}\left\{ {\zeta \left( {p + 2} \right) + \bar \zeta \left( {p + 2} \right)} \right\} - {\left( { - 1} \right)^p}{\zeta^\star}\left( {{p+2}, \overline 1} \right)  \\
&\quad - \frac{1}{2}{\ln ^2}2\left\{ {\zeta \left( {p + 1} \right) + \bar \zeta \left( {p + 1} \right)} \right\} - \ln 2{\left( { - 1} \right)^{p - 1}}\left({\zeta^\star}\left( {{p+1},  1} \right)-{\zeta^\star}\left( {\overline{p+1},  1} \right)\right)
  \\
&\quad + \sum\limits_{i = 1}^{p - 1} {{{\left( { - 1} \right)}^{i - 1}}\bar \zeta \left( {p + 1 - i} \right)\bar \zeta \left( {i + 2} \right)} - {\left( { - 1} \right)^p}\sum\limits_{n = 1}^\infty  {\frac{{{H_n}{L_n}\left( 1 \right)}}{{{n^{p + 1}}}}} .\tag{2.27}
\end{align*}
\end{thm}
\pf Using integration by parts, we can verify that
\[\int\limits_0^x {{t^{n - 1}}{\rm Li}{_q}\left( t \right)dt}  = \sum\limits_{i = 1}^{q - 1} {{{\left( { - 1} \right)}^{i - 1}}\frac{{{x^n}}}{{{n^i}}}{\rm Li}{_{q + 1 - i}}\left( x \right)}  + \frac{{{{\left( { - 1} \right)}^q}}}{{{n^q}}}\ln \left( {1 - x} \right)\left( {{x^n} - 1} \right) - \frac{{{{\left( { - 1} \right)}^q}}}{{{n^q}}}\left( {\sum\limits_{k = 1}^n {\frac{{{x^k}}}{k}} } \right),\tag{2.28}\]
\[\int_0^x {{t^{n - 1}}{{\ln }^2}\left( {1 - t} \right)} dt = \frac{1}{n} {{(x^n-1)}{{\ln }^2}\left( {1 - x} \right) } - \frac{2}{n}\sum\limits_{k = 1}^n {\frac{1}{k}} \left\{ {{x^k}\ln \left( {1 - x} \right) - \sum\limits_{j = 1}^k {\frac{{{x^j}}}{j}}  - \ln \left( {1 - x} \right)} \right\}.\tag{2.29}\]
Now, we consider the following integral
\[\int\limits_0^1 {\frac{{{\rm Li}{_p}\left( x \right){{\ln }^2}\left( {1 + x} \right)}}{x}} dx,\ 1\leq p\in Z.\]
Using (1.12) and (2.28), we can find that
\begin{align*}
&\int\limits_0^1 {\frac{{{\rm Li}{_p}\left( x \right){{\ln }^2}\left( {1 + x} \right)}}{x}} dx = 2\sum\limits_{n = 1}^\infty  {\left\{ {\frac{{{{\left( { - 1} \right)}^{n - 1}}}}{{{n^2}}} - \frac{{{H_n}}}{n}{{\left( { - 1} \right)}^{n - 1}}} \right\}\int\limits_0^1 {{x^{n - 1}}{\rm Li}{_p}\left( x \right)} dx}  \\
& =2\sum\limits_{n = 1}^\infty  {\left\{ {\frac{{{{\left( { - 1} \right)}^{n - 1}}}}{{{n^2}}} - \frac{{{H_n}}}{n}{{\left( { - 1} \right)}^{n - 1}}} \right\}\left\{ {\sum\limits_{i = 1}^{p - 1} {\frac{{{{\left( { - 1} \right)}^{i - 1}}}}{{{n^i}}}\zeta \left( {p + 1 - i} \right)}  + \frac{{{{\left( { - 1} \right)}^{p - 1}}}}{{{n^p}}}{H_n}} \right\}} \\
&=2\sum\limits_{i = 1}^{p - 1} {{{\left( { - 1} \right)}^{i - 1}}\zeta \left( {p + 1 - i} \right)\bar \zeta \left( {i + 2} \right)}  - 2{\left( { - 1} \right)^{p - 1}}{\zeta^\star}\left( {\overline{p+2},  1} \right)\\
&\quad + 2\sum\limits_{i = 1}^{p - 1} {{{\left( { - 1} \right)}^{i - 1}}\zeta \left( {p + 1 - i} \right){\zeta^\star}\left( {\overline{i+1},  1} \right) }  - 2{\left( { - 1} \right)^{p - 1}}\sum\limits_{n = 1}^\infty  {\frac{{H_n^2}}{{{n^{p + 1}}}}{{\left( { - 1} \right)}^{n - 1}}}  .\tag{2.30}
\end{align*}
On the other hand, by (2.29), we have
\begin{align*}
&\int\limits_0^1 {\frac{{{\rm Li}{_p}\left( x \right){{\ln }^2}\left( {1 + x} \right)}}{x}} dx =\sum\limits_{n = 1}^\infty  {\frac{1}{{{n^p}}}} \int\limits_0^1 {{x^{n - 1}}{{\ln }^2}\left( {1 + x} \right)} dx \\
& ={\ln ^2}2\left\{ {\zeta \left( {p + 1} \right) + \bar \zeta \left( {p + 1} \right)} \right\} - 2\ln 2\left({\zeta^\star}\left( {\overline{p+1},\overline1} \right)-{\zeta^\star}\left( {\overline{p+1},  1} \right)\right) \\
&\quad + 2{\zeta^\star}\left( {\overline{p+1},1,\overline1} \right)  .\tag{2.31}
\end{align*}
Combining (2.30) and (2.31), we can obtain (2.26). Similarly to the proof of (2.26), considering integral
\[\int\limits_0^{ - 1} {\frac{{{\rm Li}{_p}\left( x \right){{\ln }^2}\left( {1 - x} \right)}}{x}} dx,\]
we deduce (2.27) holds.
\section{Representation of Euler sums by zeta values and linear sums}
In this section, we consider the analytic representations of quadratic Euler sums of the form
\[\sum\limits_{n = 1}^\infty  {\frac{{H_n^2}}{{{n^{2m}}}}} {\left( { - 1} \right)^{n - 1}},\sum\limits_{n = 1}^\infty  {\frac{{L_n^2\left( 1 \right)}}{{{n^{2m}}}}} {\left( { - 1} \right)^{n - 1}},\sum\limits_{n = 1}^\infty  {\frac{{{H_n}{L_n}\left( 1 \right)}}{{{n^{2m}}}}} {\left( { - 1} \right)^{n - 1}},\sum\limits_{n = 1}^\infty  {\frac{{{H_n}{L_n}\left( 1 \right)}}{{{n^{2m}}}}} ,\sum\limits_{n = 1}^\infty  {\frac{{L_n^2\left( 1 \right)}}{{{n^{2m}}}}} \]through zeta values and linear sums, and give explicit formulae for several 6th-order quadratic sums in terms of zeta values and linear sums.
\begin{thm}
For $ m\geq 1$ and $m\in \Z$, we have
\begin{align*}
\sum\limits_{n = 1}^\infty  {\frac{{H_n^2}}{{{n^{2m}}}}{{\left( { - 1} \right)}^{n - 1}}}
 &=-{\zeta^\star}\left( {\overline{2m},2} \right) - \left( {2m + 1} \right){\zeta^\star}\left( {\overline{2m+1},1} \right)  - {\zeta^\star}\left( {{2m+1},1} \right)
\nonumber \\
&\quad - \zeta \left( 2 \right)\bar \zeta \left( {2m} \right) - 2\sum\limits_{l = 2}^m {{\zeta^\star}\left( {{2l-1},1} \right)\bar \zeta \left( {2m + 2 - 2l} \right)} ,\tag{3.1}
\end{align*}
\begin{align*}
\sum\limits_{n = 1}^\infty  {\frac{{{H_n}{L_n}\left( 1 \right)}}{{{n^{2m}}}}{{\left( { - 1} \right)}^{n - 1}}}
 &=\sum\limits_{l = 1}^{m - 1} {{\zeta^\star}\left( {{2l},1} \right)\bar \zeta \left( {2m + 1 - 2l} \right)}  - \ln 2{\zeta^\star}\left( {\overline{2m},1} \right)
\nonumber \\
&\quad  + \ln 2{\zeta^\star}\left( {{2m},1} \right)  + \frac{1}{2}{\zeta ^*}\left( {{2m+1},1} \right)  + {\zeta^\star}\left( {\overline{2m},2} \right)
\nonumber \\
&\quad + \frac{{2m - 1}}{2}{\zeta^\star}\left( {\overline{2m+1},1} \right)  - \frac{1}{2}\zeta \left( 2 \right)\bar \zeta \left( {2m} \right)  .\tag{3.2}
\end{align*}
\end{thm}
\pf Replacing $m$ by $2m-1$ in (2.1) and taking $s=2m$ in (2.5), we can obtain (3.1) and (3.2).
\begin{thm}
For $p\geq 1$ and $p\in \Z$, we have
\begin{align*}
&\frac{3}{2}\sum\limits_{n = 1}^\infty  {\frac{{L_n^2\left( 1 \right)}}{{{n^{p + 1}}}}} {\left( { - 1} \right)^{n - 1}} - \left( {1 + {{\left( { - 1} \right)}^{p - 1}}} \right)\sum\limits_{n = 1}^\infty  {\frac{{H_n^2}}{{{n^{p + 1}}}}{{\left( { - 1} \right)}^{n - 1}}}   \\
& =\ln 2\left({\zeta^\star}\left( {\overline{p+1},\overline1} \right)-{\zeta^\star}\left( {\overline{p+1},1} \right)\right) + \frac{1}{2}\zeta \left( 2 \right)\left\{ {\zeta \left( {p + 1} \right) + \bar \zeta \left( {p + 1} \right)} \right\}\\
&\quad + \frac{1}{2}{\zeta^\star}\left( {\overline{p+1},2} \right)  - \sum\limits_{j = 2}^p {\left( {1 - {{\left( { - 1} \right)}^{j - 1}}} \right){\zeta^\star}\left( {\overline{j},1} \right)\zeta \left( {p + 2 - j} \right)} \\
&\quad - \sum\limits_{j = 2}^{p + 1} {{\zeta^\star}\left( {{j},\overline1} \right)\bar \zeta \left( {p + 2 - j} \right)}  + \left( {p + 2 + {{\left( { - 1} \right)}^{p - 1}}} \right){\zeta^\star}\left( {\overline{p+2},1} \right)  - {\zeta^\star}\left( {{p+2},\overline1} \right)  \\
&\quad- \sum\limits_{i = 1}^{p - 1} {{{\left( { - 1} \right)}^{i - 1}}\zeta \left( {p + 1 - i} \right)\bar \zeta \left( {i + 2} \right)}  - \left( {p + 1} \right){\zeta^\star}\left( {\overline{p+2},\overline1} \right) + \ln 2{\zeta^\star}\left( {\overline{p+1},\overline1} \right).  \tag{3.3}
\end{align*}
\end{thm}
\pf From (1.8)(2.18) and (2.21), we obtain
\begin{align*}
{\zeta^\star}\left( {\overline{s},1,\overline1} \right)&=\sum\limits_{n = 1}^\infty  {\frac{{{{\left( { - 1} \right)}^{n - 1}}}}{{{n^s}}}\left( {\sum\limits_{k = 1}^n {\frac{{{L_k}\left( 1 \right)}}{k}} } \right)} \\
& =\sum\limits_{n = 1}^\infty  {\frac{{H_n^2}}{{{n^s}}}{{\left( { - 1} \right)}^{n - 1}}}  + 2\ln 2\left({\zeta^\star}\left( {\overline{s},\overline1} \right)-{\zeta^\star}\left( {\overline{s},1} \right)\right)  - \frac{3}{2}\sum\limits_{n = 1}^\infty  {\frac{{L_n^2\left( 1 \right)}}{{{n^s}}}} {\left( { - 1} \right)^{n - 1}} \\
&\quad - \sum\limits_{j = 1}^{s - 1} {{\zeta^\star}\left( {\overline{j},1} \right)\zeta \left( {s + 1 - j} \right)}  - \sum\limits_{j = 2}^s {{\zeta^\star}\left( {{j},\overline1} \right)\bar \zeta \left( {s + 1 - j} \right)}\\
&\quad+ \left( {s + 1} \right){\zeta^\star}\left( {\overline{s+1},1} \right)  - s{\zeta^\star}\left( {\overline{s+1},\overline1} \right)  + \ln 2{\zeta^\star}\left( {\overline{s},\overline1} \right)  \\
&\quad-{\zeta^\star}\left( {{s+1},\overline1} \right)+ \frac{1}{2}{\zeta^\star}\left( {\overline{s},2} \right) + \frac{1}{2}\left\{ {\zeta \left( 2 \right) - {{\ln }^2}2} \right\}\bar \zeta \left( s \right) .\tag{3.4}
\end{align*}
Taking $s=p+1$ in (3.4) and combining (2.26), we can obtain (3.3).
\begin{thm}
For $p\geq 1$ and $p \in \Z$, we have
\begin{align*}
&\sum\limits_{n = 1}^\infty  {\frac{{L_n^2\left( 1 \right)}}{{{n^{p + 1}}}}}  + \left( {{{\left( { - 1} \right)}^{p - 1}} - 1} \right)\sum\limits_{n = 1}^\infty  {\frac{{{H_n}{L_n}\left( 1 \right)}}{{{n^{p + 1}}}}} \\
& =  {\zeta^\star}\left( {{p+1},\overline2} \right) - \ln 2{\zeta^\star}\left( {{p+1},\overline1} \right)  + {\zeta^\star}\left( {\overline{p+2},\overline1} \right) + \frac{1}{2}\zeta \left( 2 \right)\zeta \left( {p + 1} \right)\\
&\quad + \sum\limits_{i = 1}^{p + 1} {{\zeta^\star}\left( {\overline{i},\overline1} \right)} \bar \zeta \left( {p + 2 - i} \right) + \left( {p + 1 + {{\left( { - 1} \right)}^p}} \right){\zeta^\star}\left( {{p+2},\overline1} \right)  \\
&\quad + \frac{1}{2}{\ln ^2}2\left\{ {2\zeta \left( {p + 1} \right) + \bar \zeta \left( {p + 1} \right)} \right\} + \ln 2{\left( { - 1} \right)^{p - 1}}\left({\zeta^\star}\left( {{p+1},1} \right)-{\zeta^\star}\left( {\overline{p+1},1} \right)\right)\\
&\quad- \sum\limits_{i = 1}^{p - 1} {{{\left( { - 1} \right)}^{i - 1}}\bar \zeta \left( {p + 1 - i} \right)} {\zeta^\star}\left( {\overline{i+1},1} \right) - \ln 2\left({\zeta^\star}\left( {{p+1},1} \right)-{\zeta^\star}\left( {{p+1},\overline1} \right)\right)\\
&\quad - \sum\limits_{i = 1}^{p - 1} {{{\left( { - 1} \right)}^{i - 1}}} \bar \zeta \left( {p + 1 - i} \right)\bar \zeta \left( {i + 2} \right) - \ln 2{\left( { - 1} \right)^{p - 1}}\left\{ {\zeta \left( {p + 2} \right) + \bar \zeta \left( {p + 2} \right)} \right\}.\tag{3.5}
\end{align*}
\end{thm}
\pf  The following identity is easily derived
\[\sum\limits_{k = 1}^n {\frac{{{H_k}}}{k}{{\left( { - 1} \right)}^{k - 1}}}  + \sum\limits_{k = 1}^n {\frac{{{L_k}\left( 1 \right)}}{k}}  = {H_n}{L_n}\left( 1 \right) + {L_n}\left( 2 \right).\tag{3.6}\]
Multiplying (3.6) by $\frac{1}{n^{p+1}}$ and summing with respect to $n$, the result is
\begin{align*}
{\zeta^\star}\left( {{p+1},\overline1,1} \right)+{\zeta^\star}\left( {{p+1},1,\overline1} \right)&=-\sum\limits_{n = 1}^\infty  {\frac{1}{{{n^{p + 1}}}}\left( {\sum\limits_{k = 1}^n {\frac{{{H_k}}}{k}{{\left( { - 1} \right)}^{k - 1}}} } \right)}  - \sum\limits_{n = 1}^\infty  {\frac{1}{{{n^{p + 1}}}}\left( {\sum\limits_{k = 1}^n {\frac{{{L_k}\left( 1 \right)}}{k}} } \right)} \\ &= -\sum\limits_{n = 1}^\infty  {\frac{{{H_n}{L_n}\left( 1 \right)}}{{{n^{p + 1}}}}}  + {\zeta^\star}\left( {{p+1},\overline2} \right).\tag{3.7}
\end{align*}
Substituting (2.19)(2.27) into (3.7), we deduce Theorem 3.3 holds.
\begin{thm}
For $l\geq 2,p\geq 1$ and $l,p\in \Z$, we have
$$\sum\limits_{n = 1}^\infty  {\frac{{\zeta _n^2\left( p \right)}}{{{n^l}}}}  = {2^{2p + l - 3}}\left\{ {\sum\limits_{n = 1}^\infty  {\frac{{{{\left( {{\zeta _n}\left( p \right) - {L_n}\left( p \right)} \right)}^2}}}{{{n^l}}}} \left( {1 - {{\left( { - 1} \right)}^{n - 1}}} \right)} \right\},\eqno(3.8)$$
$$\sum\limits_{n = 1}^\infty  {\frac{{\zeta _n^2\left( p \right)}}{{{{\left( {2n + 1} \right)}^l}}}}  = {2^{2p - 3}}\left\{ {\sum\limits_{n = 1}^\infty  {\frac{{{{\left( {{\zeta _n}\left( p \right) - {L_n}\left( p \right)} \right)}^2}}}{{{n^l}}}\left( {1 + {{\left( { - 1} \right)}^{n - 1}}} \right)} } \right\}.\eqno(3.9)$$
\end{thm}
\pf  By the definitions of harmonic numbers and alternating harmonic numbers, we have the relations
$${\zeta _{2n}}\left( p \right) - {L_{2n}}\left( p \right) = \frac{1}{{{2^{p - 1}}}}{\zeta _n}\left( p \right),{\zeta _{2n + 1}}\left( p \right) - {L_{2n + 1}}\left(p \right) = \frac{1}{{{2^{p - 1}}}}{\zeta _n}\left( p \right),\eqno(3.10)$$
$${\zeta _{2n}}\left( p \right) + {L_{2n}}\left( p \right) = 2{h_n}\left( p \right),{\zeta _{2n + 1}}\left( p \right) + {L_{2n + 1}}\left( p \right) = 2{h_{n + 1}}\left( p \right).\eqno(3.11)$$
Here ${h_n}\left( p \right)$ is defined by
$${h_n}\left( p \right) = \sum\limits_{j = 1}^n {\frac{1}{{{{\left( {2j - 1} \right)}^p}}}} .$$
Hence, from (3.10) and (3.11), we get
\[\sum\limits_{n = 1}^\infty  {\frac{{{{\left( {{\zeta _n}\left( p \right) - {L_n}\left( p \right)} \right)}^2}}}{{{n^l}}}}  = \frac{1}{{{2^{2p + l - 2}}}}\sum\limits_{n = 1}^\infty  {\frac{{\zeta _n^2\left(p \right)}}{{{n^l}}}}  + \frac{1}{{{2^{2p - 2}}}}\sum\limits_{n = 1}^\infty  {\frac{{\zeta _n^2\left( p \right)}}{{{{\left( {2n + 1} \right)}^l}}}},\tag{3.12} \]
\[\sum\limits_{n = 1}^\infty  {\frac{{{{\left( {{\zeta _n}\left( p \right) - {L_n}\left( p \right)} \right)}^2}}}{{{n^l}}}{{\left( { - 1} \right)}^{n - 1}}}  = \frac{1}{{{2^{2p - 2}}}}\sum\limits_{n = 1}^\infty  {\frac{{\zeta _n^2\left( p \right)}}{{{{\left( {2n + 1} \right)}^l}}}}  - \frac{1}{{{2^{2p + l - 2}}}}\sum\limits_{n = 1}^\infty  {\frac{{\zeta _n^2\left( p\right)}}{{{n^l}}}}.\tag{3.13} \]
Combining (3.12) and (3.13), we can obtain (3.8) and (3.9). Taking $p=1, l=2m$ in (3.8) and (3.9), we have
\[\sum\limits_{n = 1}^\infty  {\frac{{{H_n}{L_n}\left( 1 \right)}}{{{n^{2m}}}}}  = \left( {\frac{1}{2} - \frac{1}{{{2^{2m}}}}} \right)\sum\limits_{n = 1}^\infty  {\frac{{H_n^2}}{{{n^{2m}}}}}  + \frac{1}{2}\sum\limits_{n = 1}^\infty  {\frac{{L_n^2\left( 1 \right)}}{{{n^{2m}}}}}  - \frac{1}{2}\sum\limits_{n = 1}^\infty  {\frac{{H_n^2 - 2{H_n}{L_n}\left( 1 \right) + L_n^2\left( 1 \right)}}{{{n^{2m}}}}{{\left( { - 1} \right)}^{n - 1}}} ,\tag{3.14}\]
\[\sum\limits_{n = 1}^\infty  {\frac{{H_n^2}}{{{{\left( {2n + 1} \right)}^{2m}}}}}  = \frac{1}{2}\left\{ {\sum\limits_{n = 1}^\infty  {\frac{{H_n^2 - 2{H_n}{L_n}\left( 1 \right) + L_n^2\left( 1 \right)}}{{{n^{2m}}}}\left( {1 + {{\left( { - 1} \right)}^{n - 1}}} \right)} } \right\}.\tag{3.15}\]
From (2.6)(3.1)-(3.3)(3.5)(3.14) and (3.15), we have the following Theorem.
\begin{thm}
For $m \geq 1$ and $m \in Z$, the quadratic sums
\[\sum\limits_{n = 1}^\infty  {\frac{{L_n^2\left( 1 \right)}}{{{n^{2m}}}}} {\left( { - 1} \right)^{n - 1}},\sum\limits_{n = 1}^\infty  {\frac{{{H_n}{L_n}\left( 1 \right)}}{{{n^{2m}}}}} ,\sum\limits_{n = 1}^\infty  {\frac{{{H_n}{L_n}\left( 1 \right)}}{{{n^{2m}}}}} {\left( { - 1} \right)^{n - 1}}\]
\[\sum\limits_{n = 1}^\infty  {\frac{{H_n^2}}{{{n^{2m}}}}} ,\sum\limits_{n = 1}^\infty  {\frac{{H_n^2}}{{{n^{2m}}}}} {\left( { - 1} \right)^{n - 1}},\sum\limits_{n = 1}^\infty  {\frac{{L_n^2\left( 1 \right)}}{{{n^{2m}}}}} ,\sum\limits_{n = 1}^\infty  {\frac{{H_n^2}}{{{{\left( {2n + 1} \right)}^{2m}}}}} \]
are reducible to linear sums.
\end{thm}
Next, we give some explicit formulae for several 4th and 6th-order quadratic
sums in terms of zeta values and linear sums.
\begin{exa} Some illustrative examples follow.
\begin{align*}
&\sum\limits_{n = 1}^\infty  {\frac{{{H_n}{L_n}\left( 1 \right)}}{{{n^2}}}}  = \frac{{43}}{{16}}\zeta (4) + \frac{3}{4}\zeta (2){\ln ^2}2 - \frac{1}{8}{\ln ^4}2 - 3{\rm{L}}{{\rm{i}}_4}\left( {\frac{1}{2}} \right), \\
&\sum\limits_{n = 1}^\infty  {\frac{{L_n^2\left( 1 \right)}}{{{n^2}}}}  =  - \frac{{13}}{8}\zeta (4) + \frac{5}{2}\zeta (2){\ln ^2}2 + \frac{1}{{12}}{\ln ^4}2 + 2{\rm{L}}{{\rm{i}}_4}\left( {\frac{1}{2}} \right),\\
&\sum\limits_{n = 1}^\infty  {\frac{{{H_n}{L_n}\left( 1 \right)}}{{{n^2}}}{{\left( { - 1} \right)}^{n - 1}}}  = \frac{{29}}{{16}}\zeta \left( 4 \right) + \frac{3}{4}\zeta \left( 2 \right){\ln ^2}2 - 3{\rm{L}}{{\rm{i}}_4}\left( {\frac{1}{2}} \right) - \frac{1}{8}{\ln ^4}2, \\
&\sum\limits_{n = 1}^\infty  {\frac{{H_n^2}}{{{n^2}}}} {\left( { - 1} \right)^{n - 1}} = \frac{{41}}{{16}}\zeta \left( 4 \right) + \frac{1}{2}\zeta \left( 2 \right){\ln ^2}2 - \frac{1}{{12}}{\ln ^4}2 - \frac{7}{4}\zeta \left( 3 \right)\ln 2 - 2{\rm{L}}{{\rm{i}}_4}\left( {\frac{1}{2}} \right), \\
&\sum\limits_{n = 1}^\infty  {\frac{{L_n^2\left( 1 \right)}}{{{n^2}}}} {\left( { - 1} \right)^{n - 1}} =  - \frac{{41}}{{16}}\zeta \left( 4 \right) + 2\zeta \left( 2 \right){\ln ^2}2 + \frac{1}{6}{\ln ^4}2 + \frac{7}{4}\zeta \left( 3 \right)\ln 2 + 4{\rm{L}}{{\rm{i}}_4}\left( {\frac{1}{2}} \right),\\
&\sum\limits_{n = 1}^\infty  {\frac{{H_n^2}}{{{{\left( {2n + 1} \right)}^2}}}}  =  - \frac{{61}}{{16}}\zeta \left( 4 \right) + \frac{7}{2}\zeta \left( 2 \right){\ln ^2}2 + \frac{1}{3}{\ln ^4}2 + 8{\rm{L}}{{\rm{i}}_4}\left( {\frac{1}{2}} \right).
\end{align*}
\end{exa}
Using mathematica, the numerical values of 6th-order quadratic sums and 6th-order linear sums, to 30 decimal digits, are:
\begin{align*}
&\sum\limits_{n = 1}^\infty  {\frac{{{H_n}}}{{{n^5}}}{{\left( { - 1} \right)}^{n - 1}}}  =-{\zeta^\star}\left( {\bar 5,1} \right)= 0.959151942504318157165421137321\ldots\\
&\sum\limits_{n = 1}^\infty  {\frac{{{L_n}\left( 1 \right)}}{{{n^5}}}}=-{\zeta^\star}\left( { 5,\bar1} \right)  = 1.02005194570145237930331996837\ldots\\
&\sum\limits_{n = 1}^\infty  {\frac{{{L_n}\left( 1 \right)}}{{{n^5}}}{{\left( { - 1} \right)}^{n - 1}}} ={\zeta^\star}\left( { \bar5,\bar1} \right) = 0.987441426403299713771650007985\ldots\\
&\sum\limits_{n = 1}^\infty  {\frac{{{\zeta _n}\left( 2 \right)}}{{{n^4}}}{{\left( { - 1} \right)}^{n - 1}}}= -{\zeta^\star}\left( { \bar4,2} \right) = 0.934707899349253255197542851216\ldots\\
&\sum\limits_{n = 1}^\infty  {\frac{{{L_n}\left( 2 \right)}}{{{n^4}}}}=-{\zeta^\star}\left( {4,\bar2} \right)  = 1.06358224101814909880154833539\ldots\\
&\sum\limits_{n = 1}^\infty  {\frac{{H_n^2}}{{{n^4}}}{{\left( { - 1} \right)}^{n - 1}}}  = 0.889343140860204925167721149031\ldots\\
&\sum\limits_{n = 1}^\infty  {\frac{{L_n^2\left( 1 \right)}}{{{n^4}}}{{\left( { - 1} \right)}^{n - 1}}}  = 0.992525582596587194230080495889\ldots\\
&\sum\limits_{n = 1}^\infty  {\frac{{L_n^2\left( 1 \right)}}{{{n^4}}}}  = 1.02740836872025737362228676589\ldots\\
&\sum\limits_{n = 1}^\infty  {\frac{{{H_n}{L_n}\left( 1 \right)}}{{{n^4}}}}  = 1.07692404770840527955071456146\ldots\\
&\sum\limits_{n = 1}^\infty  {\frac{{{H_n}{L_n}\left( 1 \right)}}{{{n^4}}}{{\left( { - 1} \right)}^{n - 1}}}  = 0.969581748999225091836494168654\ldots\\
&\sum\limits_{n = 1}^\infty  {\frac{{{L_n}\left( 1 \right){\zeta _n}\left( 2 \right)}}{{{n^3}}}}  = 1.15935334356951415975457027807\ldots
\end{align*}
From Theorem 3.1-3.4, we find that the 6th-order quadratic sums satisfies a relation involving homogeneous combinations of $\zeta(2),\zeta(3),\zeta(4),\zeta(5),\zeta(6),\ln2$ and 6th-order linear sums.
\begin{exa}
\begin{align*}
&\sum\limits_{n = 1}^\infty  {\frac{{H_n^2}}{{{n^4}}}{{\left( { - 1} \right)}^{n - 1}}}  =  - \frac{{175}}{{32}}\zeta \left( 6 \right) + \frac{1}{2}{\zeta ^2}\left( 3 \right) - {\zeta^\star}\left( {\bar 4,2} \right)- 5{\zeta^\star}\left( {\bar 5,1} \right), \\
&\sum\limits_{n = 1}^\infty  {\frac{{L_n^2\left( 1 \right)}}{{{n^4}}}}  = \frac{{15}}{4}{\ln ^2}2\zeta \left( 4 \right) + \frac{9}{4}\zeta \left( 2 \right)\zeta \left( 3 \right)\ln 2 - \frac{{31}}{4}\zeta \left( 5 \right)\ln 2 + \frac{{231}}{{64}}\zeta \left( 6 \right) - \frac{{45}}{{32}}{\zeta ^2}\left( 3 \right)\\
& \quad \quad\quad\quad\quad\quad+ {\zeta^\star}\left( {\bar 5,\bar1} \right)  +{\zeta^\star}\left( { 4,\bar2} \right), \\
&\sum\limits_{n = 1}^\infty  {\frac{{L_n^2\left( 1 \right)}}{{{n^4}}}{{\left( { - 1} \right)}^{n - 1}}}  = \frac{{15}}{4}{\ln ^2}2\zeta \left( 4 \right) + \frac{9}{4}\zeta \left( 2 \right)\zeta \left( 3 \right)\ln 2 - \frac{{93}}{{16}}\zeta \left( 5 \right)\ln 2 + \frac{{35}}{{64}}\zeta \left( 6 \right) - \frac{{15}}{{16}}{\zeta ^2}\left( 3 \right) \\
&\quad \quad\quad\quad\quad\quad\quad \quad\quad- {\zeta^\star}\left( {\bar 4,2} \right) , \\
&\sum\limits_{n = 1}^\infty  {\frac{{{H_n}{L_n}\left( 1 \right)}}{{{n^4}}}}  =  - \frac{3}{2}\zeta \left( 2 \right)\zeta \left( 3 \right)\ln 2 + \frac{{31}}{8}\zeta \left( 5 \right)\ln 2 + \frac{{2359}}{{384}}\zeta \left( 6 \right) - \frac{7}{{64}}{\zeta ^2}\left( 3 \right) + \frac{1}{2}{\zeta^\star}\left( {\bar 5,\bar1} \right)  \\
& \quad \quad\quad\quad\quad\quad\quad+ \frac{1}{2}{\zeta^\star}\left( {4,\bar2} \right) + 2{\zeta^\star}\left( {\bar 4,2} \right) + 4{\zeta^\star}\left( {\bar 5,1} \right) ,\\
&\sum\limits_{n = 1}^\infty  {\frac{{{H_n}{L_n}\left( 1 \right)}}{{{n^4}}}{{\left( { - 1} \right)}^{n - 1}}}  =  - \frac{3}{2}\zeta \left( 2 \right)\zeta \left( 3 \right)\ln 2 + \frac{{155}}{{32}}\zeta \left( 5 \right)\ln 2 + \frac{7}{{64}}\zeta \left( 6 \right) + \frac{5}{4}{\zeta ^2}\left( 3 \right)\\
&\quad \quad\quad\quad\quad\quad\quad\quad\quad\quad\quad+ {\zeta^\star}\left( {\bar 4,2} \right)+ \frac{3}{2}{\zeta^\star}\left( {\bar 5,1} \right) .
\end{align*}
\end{exa}
On the other hand, by Cauchy product formula, we have
\[{\rm Li}_2^2\left( x \right) = 2\sum\limits_{n = 1}^\infty  {\left\{ {\frac{{{\zeta _n}\left( 2 \right)}}{{{n^2}}} + 2\frac{{{H_n}}}{{{n^3}}}} \right\}{x^n} - 6{\rm Li}{_4}\left( x \right)} .\tag{3.16}\]
Multiplying (3.16) by $\frac{{\ln \left( {1 - x} \right)}}{x}$, and integrating over $(0,x)$. The result is
\[{\rm Li}_2^3\left( x \right) = \sum\limits_{n = 1}^\infty  {\left\{ {\frac{{18}}{{{n^5}}} - 6\left[ {\frac{{{\zeta _n}\left( 2 \right)}}{{{n^3}}} + 2\frac{{{H_n}}}{{{n^4}}}} \right]} \right\}} \left\{ {{x^n}\ln \left( {1 - x} \right) - \sum\limits_{k = 1}^n {\frac{{{x^k}}}{k}}  - \ln \left( {1 - x} \right)} \right\}.\tag{3.17}\]
Taking $x=-1$ in (3.17), using Example 2, we obtain
\begin{align*}
\sum\limits_{n = 1}^\infty  {\frac{{{L_n}\left( 1 \right){\zeta _n}\left( 2 \right)}}{{{n^3}}}} &=\frac{{29}}{8}\zeta \left( 2 \right)\zeta \left( 3 \right)\ln 2 - \frac{{93}}{{32}}\zeta \left( 5 \right)\ln 2 - \frac{{1855}}{{128}}\zeta \left( 6 \right) + \frac{{17}}{{16}}{\zeta ^2}\left( 3 \right)- {\zeta^\star}\left( {\bar 5,\bar1} \right) \\
& \quad - {\zeta^\star}\left( {4,\bar2} \right) -4{\zeta^\star}\left( {\bar 4,2} \right)- 8{\zeta^\star}\left( {\bar 5,1} \right).\tag{3.18}
\end{align*}

{\bf Acknowledgments.} The authors would like to thank the anonymous
referee for his/her helpful comments, which improve the presentation
of the paper.

 \end{document}